\documentstyle[twoside]{article}
\title{An extension of a series containing Laguerre polynomials}
\author{\sc A. K. Rathie$^a$ and R. B. Paris$^b$\\
\\
${}^a\!$ Department of Mathematics, Central University of Kerala, Kasaragad 671123,\\
Kerala, India\\
E-Mail: akrathie@cukerala.edu.in\\
${}^b\!$ School of Engineering, Computing and Applied Mathematics,\\
 University of Abertay Dundee, Dundee DD1 1HG, UK\\
E-Mail: r.paris@abertay.ac.uk}
\begin{document}
\def\f#1#2{\mbox{${\textstyle \frac{#1}{#2}}$}}
\def\dfrac#1#2{\displaystyle{\frac{#1}{#2}}}
\def\boldal{\mbox{\boldmath $\alpha$}}
\newcommand{\bee}{\begin{equation}}
\newcommand{\ee}{\end{equation}}
\newcommand{\lam}{\lambda}
\newcommand{\ka}{\kappa}
\newcommand{\al}{\alpha}
\newcommand{\th}{\theta}
\newcommand{\om}{\omega}
\newcommand{\Om}{\Omega}
\newcommand{\fr}{\frac{1}{2}}
\newcommand{\fs}{\f{1}{2}}
\newcommand{\g}{\Gamma}
\newcommand{\br}{\biggr}
\newcommand{\bl}{\biggl}
\newcommand{\ra}{\rightarrow}
\renewcommand{\topfraction}{0.9}
\renewcommand{\bottomfraction}{0.9}
\renewcommand{\textfraction}{0.05}
\newcommand{\mcol}{\multicolumn}
\date{}
\maketitle
\begin{abstract}
Expressions for the summation of the series involving the Laguerre polynomials
\[S_m(\pm\nu, \pm p)\equiv e^{-x}\sum_{n=0}^\infty \frac{x^n\,L_n^{(\nu)}(x)}{(1\pm \nu\pm p)_n}\frac{(f+m)_n}{(f)_n}\]
for any non-negative integers $m$ and $p$ are obtained in terms of generalized hypergeometric functions.
These results extend previous work in the literature.
\vspace{0.4cm}

 \setcounter{section}{1}
\setcounter{equation}{0}
\vspace{0.3cm}

\end{abstract}
\vspace{0.3cm}

\begin{center}
{\bf 1. \  Introduction}
\end{center}
\setcounter{section}{1}
\setcounter{equation}{0}
\renewcommand{\theequation}{\arabic{section}.\arabic{equation}}
The generalized Laguerre polynomials $L_n^{(\nu)}(x)$ are encountered in many branches of pure and applied mathematics. They form an orthogonal set on $[0, \infty)$ with the weight function $x^\nu e^{-x}$, and can be represented as a terminating confluent hypergeometric function ${}_1F_1$ in the form
\[L_n^{(\nu)}(x)=\frac{(\nu+1)_n}{n!}\,{}_1F_1(-n; \nu+1; x),\]
where $(a)_n$ denotes the Pochhammer symbol, or rising factorial, defined by $(a)_n=\g(a+n)/\g(a)$.
The first three polynomials are given by 
\[L_0^{(\nu)}(x)=1, \quad L_1^{(\nu)}(x)=1-x+\nu, \quad  L_2^{(\nu)}(x)=\fs x^2-(\nu+2)x+\fs(\nu+1)(\nu+2).\]

In \cite{KRP1}, Kim {\it et al.\/} obtained summation formulas for the series involving the generalized Laguerre polynomial $L_n^{(\nu)}(x)$ given by
\[\sum_{n=0}^\infty \frac{x^n\,L_n^{(\nu)}(x)}{(1\pm\nu\pm p)_n}\]
for integer $p$, where $0\leq p\leq 5$. The case $p=0$ had been considered originally by Exton \cite{E}. Recently, Brychkov \cite{B} has extended these results for any integer $p\geq 0$. Alternative expressions for $S(\pm\nu, \pm p)$, also for arbitrary integer $p\geq 0$, in terms of generalized hypergeometric functions of argument $-x^2$ have been obtained in \cite{KRP2}. These results differ from, and in some cases are simpler than, those obtained in \cite{B}.
The aim of this note is to extend the above evaluations to the case 
when there is a pair of numeratorial and denominatorial parameters differing by a positive integer $m$. We derive expressions for the summation of the series given by
\[S_m(\pm\nu, \pm p)\equiv e^{-x}\sum_{n=0}^\infty \frac{x^n\,L_n^{(\nu)}(x)}{(1\pm\nu\pm p)_n}
\frac{(f+m)_n}{(f)_n}\]
for any non-negative integer $p$ and complex parameter $f\neq 0, -1, -2, \ldots$ when it is supposed that $1\pm\nu\pm p\neq 0, -1, -2, \ldots\ $.
\vspace{0.6cm}

\begin{center}
{\bf 2. \  The case $m=0$}
\end{center}
\setcounter{section}{2}
\setcounter{equation}{0}
\renewcommand{\theequation}{\arabic{section}.\arabic{equation}}
We first summarize the results obtained in \cite{KRP1} in the case $m=0$, where
\bee\label{e20}
S_0(\pm\nu,\pm p)=e^{-x}\sum_{n=0}^\infty \frac{x^n\,L_n^{(\nu)}(x)}{(1\pm\nu\pm p)_n}.
\ee
Then, for $p=0, 1, 2, \ldots\ $, we have
\bee\label{e21}
S_0(\nu,p)=\frac{(-)^p2^{2\nu+p}\g(1+\nu)}{\g(1+2\nu+p)}\sum_{s=0}^p
(-)^s\left(\!\!\begin{array}{c}p\\s\end{array}\!\!\right)\left\{A_0(\nu,p) F_0^{(1)}-\frac{4x(1+\nu)B_0(\nu,p)}{(1+\nu+p)(1+2\nu+p)}\,F_0^{(2)}\right\}
\ee
and
\bee\label{e22}
S_0(\nu,-p)=\frac{2^{2\nu-p}\g(1+\nu-p)}{\g(1+2\nu-p)}\sum_{s=0}^p
\left(\!\!\begin{array}{c}p\\s\end{array}\!\!\right)\left\{A_0(\nu,-p) F_0^{(3)}-\frac{4xB_0(\nu,-p)}{1+2\nu-p}\,F_0^{(4)}\right\},
\ee
where
\[A_0(\nu,p)=\frac{\g(\nu\!+\!\fs p\!+\!\fs s\!+\!\fs)}{\g(\fs\!-\!\fs|p|\!+\!\fs s)},\qquad B_0(\nu,p)=\frac{\g(\nu\!+\!\fs p\!+\!\fs s\!+\!1)}{\g(\fs s\!-\!\fs |p|)},\]
and
\[F_0^{(1)}={}_4F_5\left[\!\!\begin{array}{c}\vspace{0.1cm}

\fs\!+\!\fs\nu,\, 1\!+\!\fs\nu,\, \fs\!+\!\nu\!+\!\fs p\!+\!\fs s,\, \fs\!+\!\fs p\!-\!\fs s\\ \fs,\, \fs\!+\!\fs\nu\!+\!\fs p, \,1\!+\!\fs\nu\!+\!\fs p,\, \fs\!+\!\nu\!+\!\fs p,\, 1\!+\!\nu\!+\!\fs p\end{array}\!;-x^2\right],\]

\[F_0^{(2)}={}_4F_5\left[\!\!\begin{array}{c}\vspace{0.1cm}

1\!+\!\fs\nu,\, \f{3}{2}\!+\!\fs\nu,\, 1\!+\!\nu\!+\!\fs p\!+\!\fs s, \,1\!+\!\fs p\!-\!\fs s\\ \f{3}{2},\, 1\!+\!\fs\nu\!+\!\fs p,\, \f{3}{2}\!+\!\fs\nu\!+\!\fs p, 1\!+\!\nu\!+\!\fs p,\, \f{3}{2}\!+\!\nu\!+\!\fs p\end{array}\!;-x^2\right],\]

\[F_0^{(3)}={}_2F_3\left[\!\!\begin{array}{c}\vspace{0.1cm}

\fs\!+\!\nu\!-\!\fs p\!+\!\fs s,\, \fs\!+\!\fs p\!-\!\fs s\\ \fs,\, \fs\!+\!\nu\!-\!\fs p,\, 1\!+\!\nu\!-\!\fs p\end{array}\!;-x^2\right],\quad F_0^{(4)}={}_2F_3\left[\!\!\begin{array}{c}\vspace{0.1cm}

1\!+\!\nu\!-\!\fs p\!+\!\fs s,\,1\!+\!\fs p\!-\!\fs s\\ \f{3}{2},\/ 1\!+\!\nu\!-\!\fs p,\, \f{3}{2}\!+\!\nu\!-\!\fs p\end{array}\!;-x^2\right].\]
The function ${}_rF_{r+1}$ denotes the generalized hypergeometric function with $r$ numeratorial  and $r+1$ denominatorial parameters defined by
\[{}_rF_{r+1}\left[\begin{array}{c} a_1, \ldots , a_r\\b_1, \ldots , b_{r+1}\end{array}\!;z\right]
=\sum_{k=0}^\infty \frac{(a_1)_n \ldots (a_r)_n}{(b_1)_n\ldots (b_{r+1})_n}\,\frac{z^n}{n!} \qquad (|z|<\infty).\]

When the sign of the parameter $\nu$ in the denominator of (\ref{e20}) is reversed, we have for $p=0,1, 2,\ldots$
\bee\label{e23}
S_0(-\nu,p)=\frac{(-2)^p}{p!}\sum_{s=0}^p (-)^r\left(\!\!\begin{array}{c}p\\s\end{array}\!\!\right)
\left\{C_0(-\nu,p) F_0^{(5)}-\frac{4xD_0(-\nu,p)}{(1+p)(1-\nu+p)} F_0^{(6)}\right\}
\ee
and
\bee\label{e24}
S_0(-\nu,-p)=2^{-p}\sum_{s=0}^p\left(\!\!\begin{array}{c}p\\s\end{array}\!\!\right)
\left\{C_0(-\nu,-p)F_0^{(7)}-\frac{4x D_0(-\nu,-p)}{1-\nu-p} F_0^{(8)}\right\},
\ee
where
\[C_0(\nu,p)=\frac{\g(\fs\nu\!+\!\fs p\!+\!\fs s\!+\!\fs)}{\g(\fs\nu\!-\!\fs |p|\!+\!\fs s\!+\!\fs)},\qquad
D_0(\nu,p)=\frac{\g(\fs\nu\!+\!\fs p\!+\!\fs s\!+\!1)}{\g(\fs\nu\!-\!\fs |p|\!+\!\fs s)},\]
and the functions $F_0^{(j)}$ ($5\leq j\leq 8$) are defined by
\[F_0^{(5)}={}_3F_4\left[\!\!\begin{array}{c}\vspace{0.1cm}

1,\, \fs\!+\!\fs\nu\!+\!\fs p\!-\!\fs s,\, \fs\!-\!\fs \nu\!+\!\fs p\!+\!\fs s\\ \fs\!+\!\fs p,\, 1\!+\!\fs p, \,\fs\!-\!\fs\nu\!+\!\fs p,\, 1\!-\!\fs\nu\!+\!\fs p\end{array}\!;-x^2\right],\]

\[F_0^{(6)}={}_3F_4\left[\!\!\begin{array}{c}\vspace{0.1cm}

1,\, 1\!+\!\fs\nu\!+\!\fs p\!-\!\fs s,\, 1\!-\!\fs\nu\!+\!\fs p\!+\!\fs s\\1\!+\!\fs p,\, \f{3}{2}\!+\!\fs p,\, 1\!-\!\fs\nu\!+\!\fs p,\, \f{3}{2}\!-\!\fs\nu\!+\!\fs p\end{array}\!;-x^2\right],\]
 
\[F_0^{(7)}={}_2F_3\left[\!\!\begin{array}{c}\vspace{0.1cm}

\fs\!-\!\fs\nu\!-\!\fs p\!+\!\fs s,\, \fs\!+\!\fs\nu\!+\!\fs p\!-\!\fs s\\ \fs,\, \fs\!-\!\fs\nu\!-\!\fs p,\, 1\!-\!\fs\nu\!-\!\fs p\end{array}\!;-x^2\right], \]
\[F_0^{(8)}={}_2F_3\left[\!\!\begin{array}{c}\vspace{0.1cm}

 1\!-\!\fs\nu\!-\!\fs p\!+\!\fs s,\, 1\!+\!\fs\nu\!+\!\fs p\!-\!\fs s\\ \f{3}{2},\, 1\!-\!\fs\nu\!-\!\fs p,\, \f{3}{2}\!-\!\fs\nu\!-\!\fs p\end{array}\!;-x^2\right].\]
We note that $C_0(-\nu,-p)=1$ and $D_0(-\nu,-p)=\fs(s-p-\nu)$.
\vspace{0.6cm}

\begin{center}
{\bf 3. \  The case $m\neq 0$}
\end{center}
\setcounter{section}{3}
\setcounter{equation}{0}
\renewcommand{\theequation}{\arabic{section}.\arabic{equation}}
We start with the transformation \cite[(3.5)]{KRP1}
\[
e^{-x} \sum_{n=0}^\infty \frac{(a_1)_n \ldots (a_P)_n}{(b_1)_n \ldots (b_Q)_n}\,(-xy)^n\,L_n^{(\nu)}(x)=\sum_{n=0}^\infty \frac{(-x)^n}{n!}\,{}_{P+2}F_Q\left[\!\!\begin{array}{c}-n, -n-\nu, a_1, \ldots ,a_P\\b_1, \ldots , b_Q\end{array}\!;y\right],
\]
where $P$ and $Q$ are non-negative integers. 
In this, if we take $P=r$, $Q=r+1$, $y=-1$,  and $a_1=f_1+m_1, \ldots , a_r=f_r+m_r$,
$b_1=f_1, \ldots , b_r=f_r$ and $b_{r+1}=1\pm\nu\pm p$, where $m_s$ ($1\leq s\leq r$) are non-negative integers, then
\[e^{-x} \sum_{n=0}^\infty \frac{x^n\,L_n^{(\nu)}(x)}{(1\pm\nu\pm p)_n}\,\frac{(f_1+m_1)_n \ldots (f_r+m_r)_n}{(f_1)_n \ldots (f_r)_n}\hspace{4cm}\]
\bee\label{e31}
\hspace{3cm}=\sum_{n=0}^\infty \frac{(-x)^n}{n!}\,{}_{r+2}F_{r+1}\left[\!\!\begin{array}{c}-n, -n-\nu, f_1+m_1, \ldots , f_r+m_r\\1\pm\nu\pm p, f_1, \ldots , f_r\end{array}\!; -1\right].
\ee
The ${}_{r+2}F_{r+1}$ series on the right-hand side of (\ref{e31}) can be evaluated with the help of the result given in (4.2) of \cite{JCA}. However, to avoid overly complicated formulas, we restrict our attention here to the case of a single $f$-parameter ($r=1$) and adopt a different procedure. 
\vspace{0.4cm}

\noindent 3.1\ \ {\it A lemma connected with $S_m(\pm\nu,\pm p)$}
\vspace{0.2cm}

\noindent We consider the series
\bee\label{e32}
S_m(\pm\nu,\pm p)=e^{-x} \sum_{n=0}^\infty \frac{x^n\,L_n^{(\nu)}(x)}{(1\pm\nu\pm p)_n}\,\frac{(f+m)_n}{(f)_n}
=\sum_{n=0}^\infty \frac{(-x)^n}{n!}\,{}_{3}F_{2}\left[\!\!\begin{array}{c}-n, -n-\nu, f+m \\1\pm\nu\pm p, f\end{array}\!; -1\right],
\ee
for positive integer $m$ and non-negative integer $p$. Here, $x$ is an arbitrary complex variable and $f$, $\nu$
are complex parameters such that $f\neq 0, -1, -2, \ldots$ and $1\pm\nu\pm p\neq 0, -1, -2, \ldots\,$.

To proceed we make use of the following lemma that expresses the ${}_3F_2$ series in (\ref{e32}) in terms of a finite Gauss hypergeometric series.
\newtheorem{lemma}{Lemma}
\begin{lemma}
The following formula holds true
\bee\label{e32a}
S_m(\pm\nu,\pm p)=\sum_{r=0}^m \left(\!\!\begin{array}{c}m\\r\end{array}\!\!\right) \frac{x^r}{(f)_r(1\pm\nu\pm p)_r} \sum_{n=0}^\infty \frac{(-x)^n}{n!} (n+\nu+1)_r\,{}_2F_1\left[\begin{array}{c} -n, -n-\nu\\1\pm\nu\pm p+r\end{array};-1\right]
\ee
for non-negative integer $m$ and $p$.
\end{lemma}

\noindent {\it Proof.}\ \ \ We write the ${}_3F_2(-1)$ series appearing on the right-hand side of (\ref{e32})
as a finite sum of ${}_2F_1(-1)$ series by means of the result \cite[Eq.~(7.2.3.15)]{PBM}
\[{}_3F_2\left[\begin{array}{c}a, b,\\c,\end{array}\!\!\begin{array}{c}f+m\\f\end{array}\!;x\right]=
\sum_{r=0}^m \left(\!\!\begin{array}{c}m\\r\end{array}\!\!\right)\,\frac{(a)_r (b)_r}{(c)_r (f)_r}\, x^r\,
{}_2F_1\left[\begin{array}{c}a+r, b+r\\c+r\end{array}\!;x\right]\]
for positive integer $m$. Then it follows that
\bee\label{e33}
{}_{3}F_{2}\left[\!\!\begin{array}{c}-n, -n-\nu, f+m \\1\pm\nu\pm p, f\end{array}\!; -1\right]=
\sum_{r=0}^m(-)^r\left(\!\!\begin{array}{c}m\\r\end{array}\!\!\right)\frac{(-n)_r(-n-\nu)_r}{(f)_r(1\pm\nu\pm p)_r}\,{}_2F_1\left[\begin{array}{c}-n+r, -n-\nu+r\\1\pm\nu\pm p+r\end{array}\!;-1\right].
\ee
When $m>n$, we note that the sum on the right-hand side of (\ref{e33}) is over $0\leq r\leq n$, since for $n+1\leq r\leq m$ the factor $(-n)_r$ vanishes.

Substitution of (\ref{e33}) into (\ref{e32}) and use of $(-n)_r=(-)^r n(n-1)\ldots (n-r+1)$ then produces
\begin{eqnarray*}
&&\!\!\!\!\!\!\!\!S_m(\pm\nu, \pm p)\\
&=&\sum_{r=0}^m\left(\!\!\begin{array}{c}m\\r\end{array}\!\!\right)\frac{1}{(f)_r(1\pm\nu\pm p)_r}
\sum_{n=r}^\infty\frac{(-x)^n}{(n-r)!}\,(-n-\nu)_r\,{}_2F_1\left[\begin{array}{c}-n+r, -n-\nu+r\\1\pm\nu\pm p+r\end{array}\!;-1\right]\\
&=&\sum_{r=0}^m\left(\!\!\begin{array}{c}m\\r\end{array}\!\!\right)\frac{(-x)^r}{(f)_r(1\pm\nu\pm p)_r}
\sum_{n=0}^\infty\frac{(-x)^n}{n!}\,(-n-\nu-r)_r\,{}_2F_1\left[\begin{array}{c}-n, -n-\nu\\1\pm\nu\pm p+r\end{array}\!;-1\right].
\end{eqnarray*}
Finally, employing the fact that $(-n-\nu-r)_r=(-)^r(n+\nu+1)_r$ we obtain the result stated in the lemma.\ \ \ \ \ $\Box$

The ${}_2F_1$ series on the right-hand side of (\ref{e32a}) can now be evaluated with the help of the generalized Kummer summation theorem in the form \cite{RR}
\bee\label{e34a}
{}_2F_1\left[\begin{array}{c}a, b\\1\!+\!a\!-\!b\!+\!j\end{array}\!;-1\right]=2^{-2b+j}\frac{\g(b\!-\!j)\g(1\!+\!a\!-\!b\!+\!j)}{\g(b)\g(a\!-\!2b\!+\!j\!+\!1)}\sum_{s=0}^j(-)^s
\left(\!\!\begin{array}{c}j\\s\end{array}\!\!\right)\,\frac{\g(\fs a\!-\!b\!+\!\fs j\!+\!\fs s\!+\!\fs)}{\g(\fs a\!-\!\fs j\!+\!\fs s\!+\!\fs)}
\ee
and
\bee\label{e34b}
{}_2F_1\left[\begin{array}{c}a, b\\1\!+\!a\!-\!b\!-\!j\end{array}\!;-1\right]=2^{-2b-j}\frac{\g(1\!+\!a\!-\!b\!-\!j)}{\g(a\!-\!2b\!-\!j\!+\!1)}\sum_{s=0}^j\left(\!\!\begin{array}{c}j\\s\end{array}\!\!\right)\,\frac{\g(\fs a\!-\!b\!-\!\fs j\!+\!\fs s\!+\!\fs)}{\g(\fs a\!-\!\fs j\!+\!\fs s\!+\!\fs)}
\ee
for $j=0, 1, 2, \ldots\ $. If we successively put $a=-n$, $b=-n-\nu$ and $a=-n-\nu$, $b=-n$ in (\ref{e34a}) we obtain
\[{}_2F_1\left[\begin{array}{c}-n, -n\!-\!\nu\\1\!+\!\nu\!+\!j\end{array}\!;-1\right]\hspace{7cm}\]
\bee\label{e35a}
=2^{2n+2\nu+j}\frac{(-)^j \g(1\!+\!n\!+\!\nu)\g(1\!+\!\nu\!+\!j)}{\g(1\!+\!n\!+\!\nu+j)\g(n\!+\!2\nu\!+\!j\!+\!1)}\sum_{s=0}^j
(-)^s
\left(\!\!\begin{array}{c}j\\s\end{array}\!\!\right)\,\frac{\g(\fs n\!+\!\nu\!+\!\fs j\!+\!\fs s\!+\!\fs)}{\g(-\fs n\!-\!\fs j\!+\!\fs s\!+\!\fs)}
\ee
and
\bee\label{e35b}
{}_2F_1\left[\begin{array}{c}-n, -n\!-\!\nu\\1\!-\!\nu\!+\!j\end{array}\!;-1\right]=\frac{2^{2n+j}(-)^j n!}{j! (1\!+\!j)_n (1\!-\!\nu\!+\!j)_n}\sum_{s=0}^j(-)^s
\left(\!\!\begin{array}{c}j\\s\end{array}\!\!\right)\,\frac{\g(\fs n\!-\!\fs \nu\!+\!\fs j\!+\!\fs s\!+\!\fs)}{\g(-\fs n\!-\!\fs\nu\!-\!\fs j\!+\!\fs s\!+\!\fs)}
\ee
for $j=0, 1, 2, \ldots\ $. Similarly, if we successively put $a=-n$, $b=-n-\nu$ and $a=-n-\nu$, $b=-n$ in (\ref{e34b}) we obtain
\bee\label{e35c}
{}_2F_1\left[\begin{array}{c}-n, -n\!-\!\nu\\1\!+\!\nu\!-\!j\end{array}\!;-1\right]=2^{2n+2\nu-j}\frac{\g(1\!+\!\nu\!-\!j)}{\g(n\!+\!2\nu\!-\!j\!+\!1)}\sum_{s=0}^j
\left(\!\!\begin{array}{c}j\\s\end{array}\!\!\right)\,\frac{\g(\fs n\!+\!\nu\!-\!\fs j\!+\!\fs s\!+\!\fs)}{\g(-\fs n\!-\!\fs j\!+\!\fs s\!+\!\fs)}
\ee
and
\bee\label{e35d}
{}_2F_1\left[\begin{array}{c}-n, -n\!-\!\nu\\1\!-\!\nu\!-\!j\end{array}\!;-1\right]=\frac{2^{2n-j}}{(1\!-\!\nu\!-\!j)_n} \sum_{s=0}^j
\left(\!\!\begin{array}{c}j\\s\end{array}\!\!\right)\,\frac{\g(\fs n\!-\!\fs\nu\!-\!\fs j\!+\!\fs s\!+\!\fs)}{\g(-\fs n\!-\!\fs\nu\!-\!\fs j\!+\!\fs s\!+\!\fs)}
\ee
for $j=0, 1, 2, \ldots\ $.
\vspace{0.4cm}

\noindent 3.2\ \ {\it Evaluation of $S_m(\pm\nu,\pm p)$}
\vspace{0.2cm}

\noindent We consider in detail only the case corresponding to $S_m(\nu,p)$, the other cases being similar. We substitute (\ref{e35a}) into (\ref{e33})
with $j=p+r$ and use the identity $(a+n)_r=(a+r)_n (a)_r/(a)_n$ to find
\[S_m(\nu,p)=(-)^p2^{2\nu+p}\sum_{r=0}^m\left(\!\!\begin{array}{c}m\\r\end{array}\!\!\right)
\frac{(-2x)^r \g(1\!+\!\nu\!+\!r)}{(f)_r (1\!+\!\nu\!+\!p)_r \g(1\!+\!2\nu\!+\!p\!+\!r)}\sum_{s=0}^{p+r}(-)^s\left(\!\!\begin{array}{c}p\!+\!r\\s\end{array}\!\!\right)\]
\[\times\sum_{n=0}^\infty\frac{(-x)^n}{n!}\,\frac{2^{2n}(1\!+\!\nu\!+\!r)_n}{(1\!+\!\nu\!+\!p\!+\!r)_n(1\!+\!2\nu\!+\!p\!+\!r)_n}\,\frac{\g(\fs n\!+\!\nu\!+\!\fs p\!+\!\fs r\!+\!\fs s\!+\!\fs)}{\g(-\fs n\!-\!\fs p\!-\!\fs r\!+\!\fs s\!+\!\fs)}.\]
Separating even and odd powers of $x$ and making use of the identities
\[(a)_{2n}=2^{2n} (\fs a)_n(\fs a+\fs)_n,\qquad (a)_{2n+1}=2^{2n} a (\fs a+\fs)_n (\fs a+1)_n,\]
\[\frac{1}{\g(-n+a)}=\frac{(-)^n (1-a)_n}{\g(a)},\]
we find after some routine algebra the desired result\footnote{We omit the subscript $s$ on the quantities $A_r(\nu, p)$, $B_r(\nu, p)$ and $F_r^{(j)}$ to ease the notation.}
\[S_m(\nu,p)=(-)^p2^{2\nu+p}\sum_{r=0}^m\left(\!\!\begin{array}{c}m\\r\end{array}\!\!\right)
\frac{(-2x)^r \g(1\!+\!\nu\!+\!r)}{(f)_r (1\!+\!\nu\!+\!p)_r \g(1\!+\!2\nu\!+\!p\!+\!r)}\sum_{s=0}^{p+r}(-)^s\left(\!\!\begin{array}{c}p\!+\!r\\s\end{array}\!\!\right)\]
\bee\label{e36}
\times \left\{A_r(\nu,p) F_r^{(1)}-\frac{4x(1\!+\!\nu\!+\!r)B_r(\nu,p)}{(1\!+\!\nu\!+\!p\!+\!r)(1\!+\!2\nu\!+\!p\!+\!r)}\,F_r^{(2)}\right\}
\ee
for $p=0, 1, 2, \ldots\ $, where
\[A_r(\nu,p)=\frac{\g(\nu\!+\!\fs p\!+\!\fs r\!+\!\fs s\!+\!\fs)}{\g(\fs s\!-\!\fs|p|\!-\!\fs r\!+\!\fs)},\qquad
B_r(\nu,p)=\frac{\g(\nu\!+\!\fs p\!+\!\fs r\!+\!\fs s\!+\!1)}{\g(\fs s\!-\!\fs|p|\!-\!\fs r)}
\]
and
\[F_r^{(1)}={}_4F_5\left[\!\!\begin{array}{c}\vspace{0.1cm}

\fs\!+\!\fs\nu\!+\!\fs r,\, 1\!+\!\fs\nu\!+\!\fs r,\, \fs\!+\!\nu\!+\!\fs p\!+\!\fs r\!+\!\fs s,\, \fs\!+\!\fs p\!+\!\fs r\!-\!\fs s\\ \fs,\, \fs\!+\!\fs\nu\!+\!\fs p\!+\!\fs r,\, 1\!+\!\fs\nu\!+\!\fs p\!+\!\fs r,\, \fs\!+\!\nu\!+\!\fs p\!+\!\fs r,\, 1\!+\!\nu\!+\!\fs p\!+\!\fs r\end{array}\!;-x^2\right],\]

\[F_r^{(2)}={}_4F_5\left[\!\!\begin{array}{c}\vspace{0.1cm}

1\!+\!\fs\nu\!+\!\fs r,\, \f{3}{2}\!+\!\fs\nu\!+\!\fs r,\, 1\!+\!\nu\!+\!\fs p\!+\!\fs r\!+\!\fs s, \,1\!+\!\fs p\!+\!\fs r\!-\!\fs s\\ \f{3}{2},\, 1\!+\!\fs\nu\!+\!\fs p\!+\!\fs r,\, \f{3}{2}\!+\!\fs\nu\!+\!\fs p\!+\!\fs r,\, 1\!+\!\nu\!+\!\fs p\!+\!\fs r,\, \f{3}{2}\!+\!\nu\!+\!\fs p\!+\!\fs r\end{array}\!;-x^2\right].\]

An analogous procedure applied to $S_m(\nu,-p)$ produces
\[S_m(\nu,-p)=2^{2\nu-p}\g(1\!+\!\nu\!-\!p) \sum_{r=0}^m\left(\!\!\begin{array}{c}m\\r\end{array}\!\!\right)
\frac{(2x)^r(1\!+\!\nu)_r}{(f)_r \g(1\!+\!2\nu\!-\!p\!+\!r)}\sum_{s-0}^{p-r}\left(\!\!\begin{array}{c}p\!-\!r\\s\end{array}\!\!\right)\]
\bee\label{e37}
\times\left\{A_r(\nu, -p) F_r^{(3)}-\frac{4x(1\!+\!\nu\!+\!r)B_r(\nu,-p)}{(1\!+\!\nu)(1\!+\!2\nu\!-\!p\!+\!r)}F_r^{(4)}\right\}\quad (p\geq m)
\ee
for positive integer $p$ satisfying $p\geq m$, where
\[F_r^{(3)}={}_4F_5\left[\!\!\begin{array}{c}\vspace{0.1cm}

\fs\!+\!\fs\nu\!+\!\fs r,\, 1\!+\!\fs\nu\!+\!\fs r,\, \fs\!+\!\nu\!-\!\fs p\!+\!\fs r\!+\!\fs s,\, \fs\!+\!\fs p\!-\!\fs r\!-\!\fs s\\ \fs,\, \fs\!+\!\fs\nu,\, 1\!+\!\fs\nu,\, \fs\!+\!\nu\!-\!\fs p\!+\!\fs r,\, 1\!+\!\nu\!-\!\fs p\!+\!\fs r\end{array}\!;-x^2\right],\]
\[F_r^{(4)}={}_4F_5\left[\!\!\begin{array}{c}\vspace{0.1cm}

1\!+\!\fs\nu\!+\!\fs r,\, \f{3}{2}\!+\!\fs\nu\!+\!\fs r,\, 1\!+\!\nu\!-\!\fs p\!+\!\fs r\!+\!\fs s,\, 1\!+\!\fs p\!-\!\fs r\!-\!\fs s\\ \f{3}{2},\, 1\!+\!\fs\nu,\, \f{3}{2}\!+\!\fs\nu,\, 1\!+\!\nu\!-\!\fs p\!+\!\fs r,\, \f{3}{2}\!+\!\nu\!-\!\fs p\!+\!\fs r\end{array}\!;-x^2\right].\]
We observe that (\ref{e37}) is valid provided $p\geq m$. This is a consequence of (\ref{e35c}) holding for non-negative $j$, which when $j$ is set equal to $p-r$ with $0\leq r\leq m$, yields the requirement $p\geq m$. If $m>p$, then we must use (\ref{e35c}) for $0\leq r\leq p$ and (\ref{e35a}) for $p+1\leq r\leq m$; see the remark at the end of this section.

Use of (\ref{e35b}) and (\ref{e35d}) in the same manner yields the representations for $S_m(-\nu,\pm p)$ given by
\[S_m(-\nu,p)=\frac{(-2)^p}{p!} \sum_{r=0}^m\left(\!\!\begin{array}{c}m\\r\end{array}\!\!\right)
\frac{(-2x)^r(1\!+\!\nu)_r}{(f)_r (1\!-\!\nu\!+\!p)_r(1+p)_r}\sum_{s-0}^{p+r}(-)^s\left(\!\!\begin{array}{c}p\!+\!r\\s\end{array}\!\!\right)\]
\bee\label{e38}
\times\left\{C_r(-\nu, p) F_r^{(5)}-\frac{4x(1\!+\!\nu\!+\!r)D_r(-\nu,p)}{(1\!+\!\nu)(1\!-\!\nu\!+\!p\!+\!r)(1\!+\!p\!+\!r)}F_r^{(6)}\right\}
\ee
and
\[S_m(-\nu,-p)=2^{-p} \sum_{r=0}^m\left(\!\!\begin{array}{c}m\\r\end{array}\!\!\right)
\frac{(2x)^r(1\!+\!\nu)_r}{(f)_r (1\!-\!\nu\!-\!p)_r}\sum_{s-0}^{p-r}\left(\!\!\begin{array}{c}p\!-\!r\\s\end{array}\!\!\right)\]
\bee\label{e39}
\times\left\{C_{-r}(-\nu, -p) F_r^{(7)}-\frac{4x(1\!+\!\nu\!+\!r)D_{-r}(-\nu,-p)}{(1\!+\!\nu)(1\!-\!\nu\!-\!p\!+\!r)}F_r^{(8)}\right\}
\quad (p\geq m),
\ee
where
\[C_r(\nu,p)=\frac{\g(\fs\nu\!+\!\fs p\!+\!\fs |r|\!+\!\fs s\!+\!\fs)}{\g(\fs\nu\!-\!\fs|p|\!-\!\fs r\!+\!\fs s\!+\!\fs)},\quad D_r(\nu,p)=\frac{\g(\fs\nu\!+\!\fs p\!+\!\fs |r|\!+\!\fs s\!+\!1)}{\g(\fs\nu\!-\!\fs|p|\!-\!\fs r\!+\!\fs s)}.\]
Here we have defined the functions $F_r^{(j)}$ ($5\leq j\leq 8$) by
\[F_r^{(5)}={}_5F_6\left[\!\!\begin{array}{c}\vspace{0.1cm}

1,\,\fs\!+\!\fs\nu\!+\!\fs r,\, 1\!+\!\fs\nu\!+\!\fs r,\, \fs\!-\!\fs\nu\!+\!\fs p\!+\!\fs r\!+\!\fs s,\, \fs\!+\!\fs\nu\!+\!\fs p\!+\!\fs r\!-\!\fs s\\ \fs\!+\!\fs\nu,\, 1\!+\!\fs\nu,\, \fs\!+\!\fs p\!+\!\fs r,\,1\!+\!\fs p\!+\!\fs r,\,
\fs\!-\!\fs\nu\!+\!\fs p\!+\!\fs r,\, 1\!-\!\fs\nu\!+\!\fs p\!+\!\fs r\end{array}\!;-x^2\right],\]

\[F_r^{(6)}={}_5F_6\left[\!\!\begin{array}{c}\vspace{0.1cm}

1,\,1\!+\!\fs\nu\!+\!\fs r,\, \f{3}{2}\!+\!\fs\nu\!+\!\fs r,\, 1\!-\!\fs\nu\!+\!\fs p\!+\!\fs r\!+\!\fs s,\, 1\!+\!\fs\nu\!+\!\fs p\!+\!\fs r\!-\!\fs s\\ 1\!+\!\fs\nu,\, \f{3}{2}\!+\!\fs\nu, 1\!+\!\fs p\!+\!\fs r,\,\f{3}{2}\!+\!\fs p\!+\!\fs r,\,
1\!-\!\fs\nu\!+\!\fs p\!+\!\fs r,\, \f{3}{2}\!-\!\fs\nu\!+\!\fs p\!+\!\fs r\end{array}\!;-x^2\right],\]

\[F_r^{(7)}={}_4F_5\left[\!\!\begin{array}{c}\vspace{0.1cm}

\fs\!+\!\fs\nu\!+\!\fs r,\, 1\!+\!\fs\nu\!+\!\fs r,\, \fs\!-\!\fs\nu\!-\!\fs p\!+\!\fs r\!+\!\fs s,\, \fs\!+\!\fs\nu\!+\!\fs p\!-\!\fs r\!-\!\fs s\\ \fs,\, \fs\!+\!\fs\nu,\, 1\!+\!\fs\nu,\, \fs\!-\!\fs\nu\!-\!\fs p\!+\!\fs r,\, 1\!-\!\fs\nu\!-\!\fs p\!+\!\fs r\end{array}\!;-x^2\right],\]

\[F_r^{(8)}={}_4F_5\left[\!\!\begin{array}{c}\vspace{0.1cm}

1\!+\!\fs\nu\!+\!\fs r,\, \f{3}{2}\!+\!\fs\nu\!+\!\fs r,\, 1\!-\!\fs\nu\!-\!\fs p\!+\!\fs r\!+\!\fs s,\, 1\!+\!\fs\nu\!+ \!\fs p\!-\!\fs r\!-\!\fs s\\ \f{3}{2},\, 1\!+\!\fs\nu,\, \f{3}{2}\!+\!\fs\nu,\, 1\!-\!\fs\nu\!-\!\fs p\!+\!\fs r,\, \f{3}{2}\!-\!\fs\nu\!-\!\fs p\!+\!\fs r\end{array}\!;-x^2\right].\]
We note that $C_{-r}(-\nu,-p)=1$ and $D_{-r}(-\nu,-p)=\fs(s-p+r-\nu)$.

As in (\ref{e37}), we note that the representation (\ref{e39}) is valid for $p\geq m$. To deal with the situation corresponding to $m>p$ we must write the sums $S_m(\pm\nu,-p)$ as
\[S_m(\pm\nu,-p)=\sum_{r=0}^p \left(\!\!\begin{array}{c}m\\r\end{array}\!\!\right) \frac{x^r(1+\nu)_r}{(f)_r(1\pm\nu-p)_r} \sum_{n=0}^\infty \frac{(-x)^n}{n!} \frac{(1+\nu+r)_n}{(1+\nu)_n}\,{}_2F_1\left[\begin{array}{c} -n, -n-\nu\\1\pm\nu- (p-r)\end{array};-1\right]\]
\[\hspace{1.7cm}+\sum_{r=p+1}^m \left(\!\!\begin{array}{c}m\\r\end{array}\!\!\right) \frac{x^r(1+\nu)_r}{(f)_r(1\pm\nu-p)_r} \sum_{n=0}^\infty \frac{(-x)^n}{n!} \frac{(1+\nu+r)_n}{(1+\nu)_n}\,{}_2F_1\left[\begin{array}{c} -n, -n-\nu\\1\pm\nu+(r-p) \end{array};-1\right].\]
Then in the first series we employ the ${}_2F_1$ summations in (\ref{e35a}) and (\ref{e35c}) for $p-r\geq 0$,
and in the second series we employ the summations in (\ref{e35b}) and (\ref{e35d}) for $r-p>0$. We omit these
details.

\vspace{0.6cm}

\begin{center}
{\bf 4. \  Concluding remarks}
\end{center}
\setcounter{section}{4}
\setcounter{equation}{0}
\renewcommand{\theequation}{\arabic{section}.\arabic{equation}}
When $m=0$ the representations (\ref{e36})--(\ref{e39}) reduce to the corresponding cases given in (\ref{e21})--(\ref{e24}). In the cases $r=0$ and $r=1$ the hypergeometric functions appearing in $S_m(-\nu,p)$ and $S_m(\pm\nu, -p)$ contract to yield lower-order functions; this is not the case for $S_m(\nu,p)$ when $p\geq 2$.
For example, in the case of $S_m(-\nu,p)$ the ${}_5F_6$ functions contract to ${}_3F_4$ functions when $r=0$ and to ${}_4F_5$ functions when $r=1$.

In the case $p=0$, $m=1$, for example, we find that the ${}_4F_5$ functions contract to simpler ${}_0F_1$ functions to yield the result
\begin{eqnarray*}
S_1(\nu,0)&=&e^{-x}\sum_{n=0}^\infty\frac{x^n L_n^{(\nu)}(x)}{(1+\nu)_n}\,\frac{(f+1)_n}{(f)_n}\\
&=&\left(1+\frac{x}{f}\right)\,{}_0F_1\left[\begin{array}{c} -\!\!-\\1+\nu\end{array}\!;-x^2\right]-\frac{x^2}{(1+\nu)f}\,{}_0F_1\left[\begin{array}{c} -\!\!-\\2+\nu\end{array}\!;-x^2\right]\\
&=&\g(1+\nu) x^{-\nu}\left\{\left(1+\frac{x}{f}\right) J_\nu(2x)-\frac{x}{f}\,J_{\nu+1}(2x)\right\},
\end{eqnarray*}
where we have employed the duplication formula for the gamma function and expreesed  ${}_0F_1$  in terms of the $J$-Bessel function..

Finally, we mention that the summations $S_m(\pm\nu, \pm p)$ in (\ref{e36})--(\ref{e39}) have been verified numerically with the help of {\it Mathematica}.

\end{document}